\documentstyle[12pt]{article}
\begin{document}

\medskip
\begin{center}
\large

{\bf Constant mean curvature surfaces via
integrable dynamical system

\bigskip

\normalsize

B. G. Konopelchenko
\footnote
{Also: Sezione INFN di Lecce and Budker Institute of Nuclear
Physics, Novosibirsk, 630090 Russia, e-mail:konopel@le.infn.it}
and
I. A. Taimanov
\footnote
{Permanent address: Institute of
Mathematics, Novosibirsk, 630090 Russia, e-mail: taimanov@math.nsk.su}
\\

\bigskip

{\em Consortium EINSTEIN,
Dipartimento di Fisica, Universita di Lecce, 73100 Lecce, Italy}
\footnote
{European Institute for Nonlinear Studies via
Transnationally Extended Interchanges}
}
\end{center}

\bigskip

\noindent{\bf Abstract}:
It is shown that the equation which describes constant mean curvature
surfaces via the generalized Weierstrass--Enneper inducing has
Ha\-mil\-to\-nian form. Its simplest finite--dimensional reduction is  the
integrable Hamiltonian system with
two degrees of freedom. This finite-dimensional
system admits $S^1$-action and classes of $S^1$-equivalence of its
trajectories are in one-to-one correspondence with different helicoidal
constant mean curvature surfaces. Thus the interpretaion of well-known
Delaunay and do Carmo--Dajzcer surfaces via integrable
finite--dimensional Hamiltonian system is established.

\vskip1.5cm

Surfaces, interfaces, fronts, and their dynamics are key ingredients
in a number of interesting phenomena in physics. They are surface
waves, growth of crystal, deformation of membranes, propagation of
flame fronts, many problems of hydrodynamics connected with motion of
boundaries between regions of differ densities and viscosities (see,
e.g., \cite{BCC,Pelce}). Quantum field theory and statistical physics
are the important customers of surfaces too (see \cite{NPW,GPW}).

Mean curvature plays special role among the characteristics of
surfaces and their dynamics in several problems both in physics and
mathematics (see, e.g., \cite{Brakke,Yau}). Surfaces of constant mean
curvature have been studied intensively during last years (see, e.g.,
\cite{Wente,PS,K}).

In present paper we discuss a new approach for construction of constant
mean curvature surfaces. This method is based on the generalized
Weierstrass--Enneper inducing (\cite{Ken,HO,Kon}). It allows to
generate constant mean curvature surfaces via integrable dynamical
system with two degrees of freedom. The relation between the
trajectories of different types and surfaces of different types is
established.

The generalization of the Weierstrass--Enneper formulas for inducing
mi\-ni\-mal surfaces have been proposed in \cite{Ken} (see also \cite{HO})
and rediscovered in different but equivalent form in
connection with integrable nonlinear equations
in \cite{Kon}. We will use notation  and formulae from \cite{Kon}.

We start with the linear system
$$
\psi_{1z}=p\psi_2,
$$
$$
\psi_{2\bar z}=-p\psi_1,\eqno{(1)}
$$
where $p(z,\bar z)$ is a real function, $\psi_1$ and $\psi_2$ are, in
general, complex functions of the complex variable $z$, and bar denotes
the complex conjugation.  By using of the solution of (1), one
introduces  the variables
$(X^1(z,\bar z),X^2(z,\bar z),X^3(z,\bar z))$ as
follows
$$
X^1 +iX^2 = 2i \int^z_{z_0} ({\bar\psi}^2_1 dz' - {\bar\psi}^2_2d{\bar
z}'),
$$
$$
X^1-iX^2= 2i \int^z_{z_0} (\psi^2_2dz' - \psi^2_1 d{\bar
z}'),\eqno{(2)}
$$
$$
X^3=-2 \int^z_{z_0} (\psi_2{\bar \psi}_1dz' + \psi_1{\bar \psi}_2
d{\bar z}').
$$
In virtue of (1) integrals (2) do not depend on the choice of the
curve of integration.

Then one treats $z,\bar z$ as local coordinates on a surface
and $(X^1,X^2,X^3)$ as coordinates of its immersion in ${\bf R}^3$.
Formulae (2) induce a surface in ${\bf R}^3$ via the solutions of
system (1). By using of the well-known formulae, one finds the
first fundamental form
$$
{\tilde \Omega} = 4 (|\psi_1|^2+|\psi_2|^2)^2 dz d{\bar z}
\eqno{(3)}
$$
and Gaussian ($K$) and mean ($H$) curvatures
$$
K=-\frac{(\log(|\psi_1|^2+|\psi_2|^2))_{z\bar z}}
{(|\psi_1|^2+|\psi_2|^2)^2}, \ \ H=\frac{p(z,\bar z)}
{|\psi_1|^2+|\psi_2|^2}.  \eqno{(4)}
$$

This type of inducing of surfaces is the generalization of the
well-known Weierstrass--Enneper inducing of minimal surfaces. Indeed,
minimal surfaces ($H\equiv 0$) correspond to $p\equiv 0$ and in this
case formulae  (2) in terms of functions
$\psi=\frac{1}{\sqrt{2}}\psi_2$ and
$\phi=\frac{1}{\sqrt{2}}{\bar \psi}_1$
are reduced to those of Weierstrass--Enneper.

In this paper we will consider the case of constant mean curvature
surfaces. In this case $p=H(|\psi_1|^2+|\psi_2|^2)$ where $H = const$
and system (1) is reduced to the following
$$
\psi_{1t}-i\psi_{1x}=2H(|\psi_1|^2+|\psi_2|^2)\psi_2,
$$
$$
\psi_{2t}+i\psi_{2x}=-2H(|\psi_1|^2+|\psi_2|^2)\psi_1,\eqno{(5)}
$$
where $z=t+ix$.

First we note that system (5) has four obvious real integrals of
motion (independent on $t$):
$$
C_+=\int dx (\psi^2_1 + \psi^2_2 + {\bar \psi}^2_1 + {\bar \psi}^2_2),
$$
$$
C_{-}=\frac{1}{i} \int dx
(\psi^2_1 + \psi^2_2 - {\bar \psi}^2_1 - {\bar \psi}^2_2),
$$
$$
P= \int dx (\psi_{1x}{\bar \psi}_2 - {\bar \psi}_1\psi_{2x}),\eqno{(6)}
$$
$$
{\cal H} = \int dx \{ \frac{i}{2}
(\psi_{1x}{\bar \psi}_2 + {\bar \psi}_1\psi_{2x}) +
H(|\psi_1|^2+|\psi_2|^2)^2 \}.
$$

Then this system is Hamiltonian, i.e. it can be represented in the
form
$$
\psi_{1t}=\{\psi_1,{\cal H}\}, \psi_{2t} = \{\psi_2,{\cal H}\}
\eqno{(7)}
$$
where the Hamiltonian ${\cal H}$ is given by (6) and the Poisson
bracket $\{,\}$ is of the form
$$
\{F_1,F_2\}= \int dx
\{(\frac{\delta F_1}{\delta \psi_1}
\frac{\delta F_2}{\delta {\bar\psi}_2} -
\frac{\delta F_1}{\delta \psi_2}
\frac{\delta F_2}{\delta {\bar\psi}_1})-
(\frac{\delta F_2}{\delta \psi_1}
\frac{\delta F_1}{\delta {\bar\psi}_2} -
\frac{\delta F_2}{\delta \psi_2}
\frac{\delta F_1}{\delta {\bar\psi}_1})
\}.\eqno{(8)}.
$$
The corresponding symplectic form is
$$
\Omega = d\psi_1 \wedge d{\bar\psi}_2  + d{\bar\psi}_1 \wedge d\psi_2
$$
and the Lagrangian is given by the following formula
$$
{\cal L} = \psi_1{\bar\psi}_{2z}-{\bar\psi}_{1{\bar z}}\psi_2 +
\frac{H}{2}(|\psi_1|^2+|\psi_2|^2)^2.
$$

Thus formula (2) establishes the correspondence between the
trajectories of the infinite--dimensional Hamiltonian system (5) and
surfaces of constant mean curvatures.

Let us put
$$
H \neq 0
$$
to omit the discussion of minimal surfaces.

Let us also restrict ourselves to the particular case of this inducing
with $p=p(t)$.  It is not difficult to show that under this
constraint the only admissible solutions, of system (5), which are
representable by finite sums of terms of the type $f(t)\exp{i\rho x}$
are of the form
$$
\psi_1 = r(t)\exp{(i\lambda x)},
\psi_2 = s(t)\exp{(i\lambda x)},
\eqno{(9)}
$$
where $\lambda (\neq 0)$ is real parameter and $r(t)=p_1+ip_2$ and
$s(t)=q_1+iq_2$ are
complex-valued functions.
System (5) in these variables has the following form
$$
r_t+\lambda r - 2H(|r|^2+|s|^2)s=0,
$$
$$
s_t-\lambda s + 2H(|r|^2+|s|^2)r=0, \eqno{(10)}
$$
or equivalent system of four equations
in terms of real and imaginary parts of $r$ and $s$.
It has the Hamiltonian form
$$
\frac{\partial p_i}{\partial t} = \{p_i,{\cal H}_0\}_0,
\frac{\partial q_j}{\partial t} = \{q_j,{\cal H}_0\}_0,\ \ i,j=1,2,
$$
with the Hamiltonian function
$$
{\cal H}_0 = \frac{H}{2}(p^2_1+p^2_2+q^2_1+q^2_2)^2-
\lambda(p_1q_1+p_2q_2)
$$
and with respect to the usual Poisson
brackets $\{,\}_0$ generated by the simplectic form
$$
\Omega_0=dp_1 \wedge dq_1 + dp_2 \wedge dq_2.
$$

It is easy to notice that the Hamiltonian function ${\cal H}_0$ can be
obtained from ${\cal H}$ by using of the finite dimensional reduction
(9).  Hamiltonian system (10) has another first integral
$$
M = p_1q_2 - p_2q_1
$$
which is in involution with the Hamiltonian ${\cal H}_0$ and moreover
these first integrals are functionally independent everywhere except
the zero ($p_i=q_j=0$). Thus we conclude that system (10) is
integrable.

This system is not only integrable but also $S^1$-symmetric. Its
Hamiltonian, the additional first integral $M$ and the Poisson
structure are preserved by the following $S^1$-action:
$$
\cases{p_1 \rightarrow p_1 \cos{\phi} - p_2\sin{\phi} &  \cr
p_2 \rightarrow p_1 \sin{\phi} + p_2 \cos{\phi} & },\ \ \ \
\cases{q_1 \rightarrow q_1 \cos{\phi} - q_2\sin{\phi} &  \cr
q_2 \rightarrow q_1 \sin{\phi} + q_2 \cos{\phi} & }.
\eqno{(11)}
$$

Let us assume without loss of generality that
$$
\lambda = H = \frac{1}{2}.
$$
Formulae (2) obtain the following form
$$
X^1=-2\int\{[(p^2_1+q^2_1-p^2_2-q^2_2)\cos{x} -
2(p_1p_2+q_1q_2)\sin{x}] dx
$$
$$
+[2(q_1q_2-p_1p_2)\cos{x}+
(q^2_1+p^2_2-q^2_2-p^2_1)\sin{x}] dt \},
$$
$$
X^2=2\int\{[2(p_1p_2+q_1q_2)\cos{x} + (p^2_1+q^2_1-p^2_2-q^2_2)\sin{x}]
dx
$$
$$+ [(p^2_1+q^2_2-p^2_2-q^2_1)\cos{x} + 2(q_1q_2-p_1p_2)\sin{x}] dt \},
\eqno{(12)}
$$
$$
X^3=-4\int\{(p_1q_1+p_2q_2) dt - (p_1q_2-p_2q_1)\} dx.
$$

Trajectories of Hamiltonian system (10) which are different modulo
symmetry (11) describe different constant mean curvature surfaces by
using of formulas (12).
It also follows from (12) that
these surfaces are invariant under the following helicoidal
transform:
$$
\cases{X^1 \rightarrow X^1 \cos{\tau} - X^2 \sin{\tau} & \cr
X^2 \rightarrow X^1 \sin{\tau} + X^2 \cos{\tau} & \cr
X^3 \rightarrow X^3 + 4M\tau} ,\eqno{(13)}
$$
and the restriction, of this transform, to the surface coincides with
the shift of $Im z = x :\ \ x \rightarrow x + \tau .$

We see that if $M=0$ then we obtain a surface of
revolution. All these surfaces are equivalent modulo (11) to surfaces
with $p_2 \equiv q_2 \equiv 0$.
It is not complicated to give a qualitative analysis of the behaviour
of the restriction of (10) onto this plane. This vector field has three
zeros at points  $(0,0)$ and $\pm\frac{1}{2},
\pm\frac{1}{2})$. The second ones correspond to
cylinders of revolution. At these points the Hamiltonian ${\cal H}_0$
is equal to $-\frac{1}{32}$. These points are bounded by cycles on which
Hamiltonian is negative but more than $-\frac{1}{32}$ and which
correspond to unduloids (i.e., the Delaunay surfaces which are embedded
into ${\bf R}^3$ and differ from cylinder and round sphere).
Hamiltonian vanishes at the zero point and two separatrices which come
from $(0,0)$ and arrives to it. These separatrices correspond to round
sphere with a pair of truncated points and bound a domain where
Hamiltonian is negative. The domain ${\cal H}_0 > 0$ is fibered by
cycles of Hamiltonian system (10) and these cycles correspond to
nodoids (i.e., Delaunay surfaces which have selfintersections).

Thus we obtain very natural Hamiltonian interpretation for the
well-known family of Delaunay surfaces (\cite{D}).

In the same manner it is shown that the full family of
surfaces which corresponds to solutions of (10) with $ M \neq 0$
coincides with the
family of helicoidal surfaces of constant mean curvature which were
constructed in \cite{CD}.

\bigskip

This work was supported by INTAS (an international association for
the promotion of cooperation with scientists from independent states of
the former Soviet Union). The second author (I.A.T.) also acknowledges
partial support by the Russian Foundation for Fundamental Studies
(grant 94-01-00528).

\end{document}